\documentclass[11pt]{amsart}
\usepackage{graphicx}
\usepackage{amsmath}
\usepackage{color}

\newtheorem{thm}{Theorem}[section]

\newtheorem{prop}[thm]{Proposition}

\theoremstyle{definition}
\newtheorem{defn}[thm]{Definition}

\theoremstyle{remark}
\newtheorem{note}[thm]{Remark}
\newtheorem{example}[thm]{Example}

\numberwithin{equation}{section}

\newcommand{\w}[1]{\widetilde{#1}}
\newcommand{\fra}[2]{\displaystyle{\frac{#1}{#2}}}
\newcommand{\R}{\mathbb{R}}
\usepackage[thinlines]{easytable}

\begin{document}

\title[Slant Submanifolds of Norden Manifolds]{Slant Submanifolds of Norden Manifolds}

\author[P. Alegre]{Pablo Alegre}
 \address{Departamento de Econom\'{\i}a, M\'etodos Cuantitativos e Historia Econ\'omica, \'Area de Estad\'{\i}stica e Investigaci\'on Operativa.  Universidad Pablo de Olavide. Ctra. de Utrera km. 1, 41013 Sevilla, Spain}

\email[Corresponding author]{psalerue@upo.es}
\author[A. Carriazo]{Alfonso Carriazo}
 \address{Departamento de Geometr\'{i}a y Topolog\'{i}a. Universidad de Sevilla. c/ Tarfia s/n, 41012 Sevilla, Spain}

\email{carriazo@us.es}

\author[F. Planas]{Francisco Planas}
 \address{Departamento de Ciencias de la Educaci\'on. Universidad Fernando III. Sevilla, Spain}

\email{francisco.planas@ceu.es}
\thanks{Both authors are partially supported by the MINECO-FEDER grant MTM2014-52197-P. They are members of the PAIDI group FQM-327 (Junta de Andaluc\'ia, Spain). The second author is also a member of the Instituto de Matem\'aticas de la Universidad de Sevilla (IMUS)}

\begin{abstract}
In this paper we introduce the notion of  slant submanifolds of a Norden manifold. We study their first properties and present a whole gallery of examples.
\end{abstract}



\maketitle

\section{Introduction}

In \cite{slantchen}, B.-Y. Chen introduced slant submanifolds of an almost Hermitian manifold, as those submanifolds for which the angle $\theta$ between $JX$ and the tangent space is constant, for any tangent vector field $X$. They play an intermediate role between complex submanifolds ($\theta=0$) and totally real ones ($\theta=\pi/2$). Since then, the study of slant submanifolds has produced an incredible amount of results and examples. Moreover, some generalizations of them have also been defined, such as semi-slant, bi-slant or generic submanifolds. 

Also the study of slant submanifolds has been extended both to odd dimensional and semi-Riemannian environments.

Consider a $(M,J,g)$ with $J$ a $(1,1)$ tensor field such as $J^2=\pm I$ and $g$ a metric, it could be considered diffrenet cases depending on the compatiblity relatioship existing between both $J$ and $g$: 

\begin{center}
\begin{tabular}{|l|r|r|}
\hline\\[-1em]
&$g(JX,JY)=g(X,Y)$&$g(JX,JY)=-g(X,Y)$\\[2pt]
\hline\\[-1em]
$J^2=-I$ & almost complex manifold & Norden manifold\\[2pt]
\hline\\[-1em]
$J^2=I$ & almost product manifold & almost para-complex manifold\\[2pt]
\hline
\end{tabular}
\end{center}



Norden manifolds were firstly called $B$-manifolds by A.P. Norden \cite{norden} and sometimes the are also called anti-Kaehler \cite{borowiec}.

Slant submanifolds of an almost product manifold were introduced by B. Sakin \cite{sahin}, the first two authors studied curves as slant submanifolds, \cite{acslantproduct}, and innitiated the study of slant submanifolds of a almost para-complex manifold, \cite{acslantpara}. The objetive of the present paper is fullfilled this panorama defining slant submanifolds of a Norden manifold. the techniques we will used are similar to the ones used in these three cited studies, but there are also some differences and similarities in the Norden case that we will highlight along the paper. 

Let us notice that the paper \cite{gupta} was a previous attempt to define slant submanifolds of a Kaehler-Norden Hermitian manifold. However, we think that it is not completely right. Actually, the definition given in that paper is not correct because of some special difficulties dealing with semi-Riemannian geometry. We will give more details in Remark \cite{gupta}. Slant submanifolds of almost poly-Norden Riemannian manifolds have also been studied, \cite{ay}, but they do not match in the case under study.

The structure of the paper is the following: in a preliminary section, the main background about Norden manifolds and submanifolds is remembered, in the following slant submanifolds of a Norden manifold are defined and characterized. In the third section some non trivial examples are given. Curves which could be considered slant submanifold in some particular Norden manifolds are studied in the fourth section.

\section{Preliminaries}

Let $(\w M,J)$ be an almost complex manifold, $dim(\w M)=2n$. A metric $g$ on $\w M$ is said to be {\em Norden} if the complex structure $J$ is an antiisometry of the tangent space, that is
\begin{equation}\label{cypc}
\begin{array}{ccc}
J^2 X= -X, & g(J X,JY)=-g(X,Y),
\end{array}
\end{equation}
for any vector fields $X,Y$ on $\widetilde M$, it is called a {\it almost complex manifold with Norden metric (Norden manifold)} and {\it almost complex manifold with B-metric}, \cite{borowiec,ganchev,ganchev2,chr}.
It is said to be {\it Kaehler Norden} or {\it anti-Kaehler} if, in addition, $\widetilde \nabla J=0$.
From (\ref{cypc}) it is deduced
\begin{equation}\label{cambio}
g(JX,Y)=g(X,JY).
\end{equation}

Let now $M$ be a submanifold of $(\w M,J,g)$. The Gauss and Weingarten formulas are given by
\begin{equation}\label{gauss}
\widetilde\nabla_XY=\nabla_XY+h(X,Y),
\end{equation}
\begin{equation}\label{Weingarten}
\widetilde\nabla_X V=-A_VX+\nabla^\perp_XV,
\end{equation}
for any tangent vector fields $X,Y$ and any normal vector field $V$, where $h$ is the second fundamental form of $M$, $A_V$ is the Weingarten endomorphism associated with $V$ and $\nabla^\perp$ is the normal connection.

%

For every tangent vector field $X$, we write
\begin{equation}\label{decomp}
J X=PX+FX,
\end{equation}
where $PX$ is the tangential component of $JX$ and $FX$ is the normal one. And for every normal vector field $V$,
$$J V=tV+fV,$$
where $tV$ and $fV$ are the tangential and normal components of $JV$, respectively.

\section{Definition and characterization results}

For an almost Hermitian manifold $(\w M,J,g)$ it holds the Cauchy-Schwarz inequality 
\begin{equation*}\label{1}
|g_p(u,v)|\leq \|u\|\|v\|,
\end{equation*}
for any tangent vector fields $u,v$ at $p$. Considering $M$ a submanifold isometrically immersed in $\w M$, we have
$$|g(JX,PX)|\leq \|JX\|\|PX\|,$$
for any tangent vector field $X$. And for any non-null tangent vector field the Wirtinger angle is given by 
$$\cos\theta=\frac{g(JX,PX)}{\|JX\|\|PX\|}=\frac{\|PX\|}{\|JX\|},$$ for certain function $\theta$. $M$ is said to be \emph{slant} if that function $\theta$ is constant, \cite{slantchen}, and then $\theta$ is called the \emph{slant angle}. 

This reasoning does not hold for light-like vector vector fields in the semi-Riemannian setting. As it was done in \cite{acslantpara} for almost product manifolds, we define slant submanifolds in para Hermitian manifolds in the following way:
\begin{defn}
A submanifold $M$ of a Norden manifold $(\w M,J,g)$ is called \emph{slant} if for every space-like or time-like tangent vector field $X$, the quotient $g(PX,PX)/g(JX,JX)$ is constant.
\end{defn}

\begin{note}
It is clear that, if $M$ is a complex submanifold, then $P\equiv J$, and so the above quotient is equal to $1$. On the other hand, if $M$ is totally real, then $P\equiv 0$ and the quotient equals $0$. Therefore, both complex and totally real submanifolds are particular cases of slant submanifolds. A neither complex nor totally real slant submanifold will be called \emph{proper slant}.
\end{note}

\

Consider $X$ a unitary spacelike vector field in $M$. 

\begin{itemize}
\item[1)] If $PX$ is time-like and $g(PX,PX)>1$, then $g(FX,FX)=-1+g(PX,PX)>0$ and so $FX$ is  space-like.
\item[2)] If $PX$ is time-like and $g(PX,PX)<1$, then $g(FX,FX)=-1+g(PX,PX)<0$ and so $FX$ is time-like.
\item[3)] Finally, if $PX$ is a space-like, $g(Fe_1,Fe_1)=-1-g(PX,PX)$, and $FX$ is a time-like vector field.
\end{itemize}

Therefore, we can distinguish three different cases:

\begin{defn}
Let $M$ be a proper slant submanifold of a Norden manifold $(\w M,J,g)$. We say that it is of
\begin{itemize}
\item[{\rm type 1}] if for any space-like (time-like) vector field $X$, $PX$ is time-like (space-like), and $\fra{|PX|}{|JX|}>1$,
\item[{\rm type 2}] if for any space-like (time-like) vector field $X$, $PX$ is time-like (space-like), and $\fra{|PX|}{|JX|}<1$,
\item[{\rm type 3}] if for any space-like (time-like) vector field $X$, $PX$ is space-like (time-like).
\end{itemize}
%
%
\end{defn}

\begin{note} \label{error}
As we have showed, we must be very careful by studying separately each case, since the Wirtinger angle can not be defined for all vectors fields. This is the main problem of \cite{gupta}, where it is defined as $\cos\Theta(X)=|PX|/|JX|$ for any tangent vector field $X$.
\end{note}

\begin{note}
The case $g(PX,PX)=0$ ($g(PX,PX)=g(JX,JX)$, respectively) for any space-like or time-like $X$,  corresponds to the totally real case (complex case), that is $P\equiv 0$ ($F\equiv 0$ or equivalently $P\equiv J$). Therefore, the case $PX$ being light-like is excluded. Indeed, let us suppose that $g(PX,PX)=0$, for any space-like or time-like tangent vector field $X$. Since any light-like vector field can be approximated by a sequence of space-like or time-like vector fields, it is clear that $g(PX,PX)=0$, for any vector field $X$, independently of its causal character. Then,
$$0=g(P(X+Y),P(X+Y))=g(PX+PY,PX+PY)=2g(PX,PY),$$
for any $X,Y$.  
Thus, from the previous equation it holds $g(PX,PY)=0$ for any tangent vector field $Y$. But $g(PX,Z)=0$ for any $Z\in D$
, where we have written $TM=P(TM)\oplus D$. Hence, $PX\in\mbox{rad}(TM,g)$ and then, since $g$ is non degenerate, $PX=0$. Therefore, $M$ is totally real. The complex case can be similarly proved.


\end{note}

\

Most of the results we obtained for slant submanifolds of a para Hermitian manifold at \cite{acslantpara} can be almost directly translated to slant submanifolds of a Norden manifold. We detail the first characterization theorem and just list the other results as the proofs are quite similar and can be omitted.
\begin{thm}\label{aannaa}
Let $M$ be a submanifold of a Norden manifold $(\w M, J,g)$. Then,
\begin{itemize}

\item[1)] $M$ is slant of type 1 if and only if for any space-like (time-like) vector field $X$, $PX$ is time-like (space-like), and there exists a constant $\lambda\in(1,+\infty)$ such that
\begin{equation}\label{anna2}
P^2=\lambda Id.
\end{equation}
We write $\lambda=-\cosh^2\theta$, with $\theta >0$.
\item[2)] $M$ is slant of type 2 if and only if for any space-like (time-like) vector field $X$, $PX$ is time-like (space-like), and there exists a constant $\lambda\in(0,1)$ such that
\begin{equation}\label{anna3}
P^2=\lambda Id.
\end{equation}
We write $\lambda=-\cos^2\theta$, with $0<\theta<2\pi$.
\item[3)] $M$ is slant of type 3 if and only if for any space-like (time-like) vector field $X$, $PX$ is space-like (time-like), and there exists a constant $\lambda\in(-\infty,0)$ such that
\begin{equation}\label{anna}
P^2=\lambda Id.
\end{equation}
We write $\lambda=\sinh^2\theta$, with $\theta >0$.
\end{itemize}
In every case, we call $\theta$ the {\rm slant angle}.
\end{thm}
\begin{proof}
In the first case, if $M$ is slant of type 1, for any space-like tangent vector field $X$, $PX$ is time-like, and, by virtue of \eqref{cypc}, $JX$ also is. Moreover, they satisfy $|PX|/|JX|>1$. So, there exists $\theta>0$ such that
\begin{equation}\label{papana}
\cosh\theta=\fra{|PX|}{|JX|}=\fra{\sqrt{-g(PX,PX)}}{\sqrt{-g(JX,JX)}}.
\end{equation}

If we now consider $PX$, then, in a similar way, we obtain:
\begin{equation}\label{mericano}
\cosh\theta=\fra{|P^2X|}{|JPX|}=\fra{|P^2X|}{|PX|}.
\end{equation}

On the one hand,
\begin{equation}\label{papana1}
g(P^2X,X)=g(JPX,X)=g(PX,JX)=g(PX,PX)=|PX|^2.
\end{equation}
Therefore, using (\ref{papana}), (\ref{mericano}) and (\ref{papana1}) $$g(P^2X,X)=|PX|^2=|P^2X||JX|=|P^2 X||X|.$$

On the other hand, 
since both $X$ and $P^2X$ are space-like, it follows that they are collinear, that is $P^2X=\lambda X$. Finally, from (\ref{papana}) 
$$\cosh^2\theta=\fra{-g(PX,PX)}{-g(JX,JX)}=\fra{-\lambda g(X,X)}{g(X,X)}=-\lambda.$$ 

Everything works in a similar way for any time-like tangent vector field $Y$, but now, $PY$ and $JY$ are space-like and so, instead of (\ref{papana}) we should write:
$$\cosh\theta=\fra{|PY|}{|JY|}=\fra{\sqrt{g(PY,PY)}}{\sqrt{g(JY,JY)}}.$$

Since $P^2X=\lambda X$, for any space-like or time-like $X$, it also holds for light-like vector fields and so we have that $P^2=\lambda Id$.

The converse is just a simple computation.

In the second case, if $M$ is slant of type 2, for any space-like or time-like vector field $X$, $|PX|/|JX|<1$, and so there exists $\theta>0$ such that
\begin{equation*}\label{papana2}
\cos\theta=\fra{|PX|}{|JX|}=\fra{\sqrt{-g(PX,PX)}}{\sqrt{-g(JX,JX)}}.
\end{equation*}

Proceeding as before, we can prove that $g(P^2X, X)=| P^2X||X|$ and, as both $X$ and $P^2X$ are space-like vector fields, it follows that they are collinear, that is $P^2X=\lambda X$. Again, the converse is a direct computation.

Finally, if $M$ is slant of type 3, for any space-like vector field $X$, $PX$ is also space-like, and there exists $\theta>0$ such that
\begin{equation*}\label{papana3}
\sinh\theta=\fra{|PX|}{|JX|}=\fra{\sqrt{g(PX,PX)}}{\sqrt{-g(JX,JX)}}.
\end{equation*}

Once more, we can prove that $g(P^2X, X)=| P^2X||X|$ and $P^2X=\lambda X$. And again, the converse is a direct computation.
\end{proof}

\begin{note}
In the second case, for slant submanifolds of type 2, the slant angle really coincides with the Wirtinger angle, i.e., the angle between $JX$ and $PX$.

One more consideration: as every light-like vector field can be decomposed as a sum of one space-like and one time-like vector field, we directly obtain that conditions (\ref{anna2}), (\ref{anna3}) and (\ref{anna}) also hold for every light-like vector field.

Finally, let us point out that for both slant submanifolds of type 1 and 2, if $X$ is space-like, then $PX$ is time-like. Therefore every  slant submanifold of type 1 or 2 must be a neutral semi-Riemannian manifold.
\end{note}

In fact, for types 1 and 2 it is only necessary to ask space-like vector fields to satisfy conditions (\ref{anna2}) and (\ref{anna3}):

\begin{prop}
Let $M$ be a submanifold of a Norden manifold $(\w M, J,g)$. Then $M$ is a slant submanifold of
\begin{itemize}
\item[{\rm type 1}] if and only if $P^2X=-\cosh^2\theta X$ for every space-like vector field $X$.
\item[{\rm type 2}]  if and only if $P^2X=-\cos^2\theta X$ for every space-like vector field $X$.
\end{itemize}
\end{prop}

\begin{proof}

The proof is similar to Theorem 3.8 at \cite{acslantpara} and hence it could be ommitted.
\end{proof}

\

The case of slant submanifolds of type 3 is completly different as they are not always neutral submanifolds, in fact they can be space-like or time-like submanifolds.

\

Now we give another characterization result:

\begin{prop}
Let $M$ be a submanifold of a Norden manifold $(\w M, J,g)$. Then,
%
%
Let $M$ be a submanifold of a para Hermitian manifold $\w M$. Then $M$ is a slant submanifold of
\begin{itemize}
\item[{\rm type 1}] if and only if $tFX=\sinh^2\theta X$ for every space-like (time-like) vector field $X$.
\item[{\rm type 2}] if and only if $tFX=-\sin^2\theta X$ for every space-like (time-like) vector field $X$.
\item[{\rm type 3}] if and only if $tFX=-\cosh^2\theta X$ for every space-like (time-like) vector field $X$.
\end{itemize}
\end{prop}

\begin{proof}

The proof is similar to Proposition 3.9 at \cite{acslantpara} and hence it could be ommitted.
\end{proof}

\


Moreover we can characterize the special case of neutral slant submanifolds of type 1 and 2, with half the dimension of the ambient space, by studying the normal fields:

\begin{prop}
Let $M_s^{2s}$ be a submanifold of a Norden manifold $(\w M^{4s}_{2s},J,g)$. Then $M$ is a slant submanifold of
\begin{itemize}
\item[{\rm type 1}] if and only if $f^2V=-\cosh^2\theta V$ for every space-like (time-like) normal vector field $V$.
\item[{\rm type 2}] if and only if $f^2V=-\cos^2\theta V$ for every space-like (time-like) normal vector field $V$.
\end{itemize}
\end{prop}

\begin{proof}

The proof is similar to Theorem 3.10 at \cite{acslantpara} and hence it could be ommitted.
\end{proof}

\

As we mentioned above, a completely different study is necessary for type 3 slant submanifolds as they are not always neutral:

\begin{thm}\label{af}
Let $M_{2i}^{2s}$ $(0<i<s)$ be a submanifold of a Norden manifold $(\w M^{4s}_{2s},J,g)$. Then, $M$ is a slant submanifold of type 3 if and only if $f^2 V=\sinh^2 \theta V$ for every normal vector field $V$.
\end{thm}

\begin{proof}

The proof is similar to Theorem 3.11 at \cite{acslantpara} and hence it could be ommitted.
\end{proof}

\section{Examples}

In this section we present examples of each type of slant submanifolds. Of course, it will be necessary to exclude the complex and totally real cases, by choosing in any case the appropriate parameters.

Let us consider $\R^4$ with the following para Kaehler structure:
$$J=\left(\begin{array}{cccc}
0 & 1 & 0 & 0\\
1 & 0 & 0 & 0\\
0 & 0 & 0 & 1\\
0 & 0 & 1 & 0
\end{array}\right),\quad\quad g=\left(\begin{array}{cccc}
1 & 0 & 0 & 0\\
0& -1 & 0 & 0\\
0 & 0&1 & 0\\
0 & 0 & 0 &-1
\end{array}\right).$$



\

The first examples we can find are of type 2:

\begin{example}
For any $a, b \in \R$,
$$x(u,v)=(u\sin a,v \sin b ,u\cos a, v\cos b)$$
defines a slant submanifold of type 2 in $(\R^4,J,g)$, with $P^2=\cos^2(a-b)Id$.
\end{example}


\begin{example}\label{cuatro}
For any $a, b \in \R$ with ,
$$x(u,v)=(u\sinh a,v \sinh b,u\cosh a, v\cosh b)$$
defines a slant submanifold of type 2 in $(\R^4,J,g)$, with $P^2=\fra{\cosh^2(a+b)}{\cosh (2a)\cosh(2b)} Id$.
\end{example}

%
%
%
%

\begin{example}
For any $a,b\in \R$ with $a^2+b^2\neq 0$,
$$x(u,v)=(au,v,bu, v)$$
defines a slant submanifold of type 2 in $(\R^4,J,g)$, with $P^2=\fra{(a+b)^2}{2(a^2+b^2)}Id$.
\end{example}

In fact we can present easy examples of all three types of slant submanifolds:

\begin{example}\label{moonlight}
For any $a,b \in \R$ with $a^2+b^2\neq 1$,
$$x(u,v)=(au, v,bu, u)$$ defines a slant submanifold in $(\R^4,J,g)$, with $P^2=\fra{a^2}{-1+a^2+b^2}Id$, and it is
\begin{itemize}
\item[1)] of type 1 if $a^2+b^2>1$ and $b^2<1$,
\item[2)] of type 2 if $a^2+b^2>1$ and $b^2>1$,
\item[3)] time-like of type 3 if $a^2+b^2<1$.
\end{itemize}
\end{example}

\

Therefore, with a proper election of $a$ and $b$ we obtain the following three nice examples:

\begin{example}
For any $\theta>0$,
$$x(u,v)=(u\cosh^2 \theta,v,u\sqrt{1-\sinh^2\theta}, u)$$
defines a slant submanifold of type 1
in $(\R^4,J,g)$, with slant angle $\theta$.
\end{example}

\begin{example}
For any $0<\theta<\pi/2$,
$$x(u,v)=(u\cos^2 \theta,v,u\sqrt{\sin^2\theta+1}, u)$$
defines a slant submanifold of type 2
in $(\R^4,J,g)$, with slant angle $\theta$.
\end{example}

\begin{example}\label{lalaland}
For any $\theta>0$,
$$x(u,v)=(u\sinh^2 \theta,v,u\sqrt{1-\cosh^2\theta}, u)$$
defines a slant submanifold of type 3
in $(\R^4,J,g)$, with slant angle $\theta$.
\end{example}


%
%
%
%

\

We can also obtain new examples based on Examples 8.4 and 8.5 from \cite{slantchen}. Now we consider $\R^4$ with a different para Kaehler structure:
$$J_1=\left(\begin{array}{cccc}
0 & 0 & 1 & 0\\
0 & 0 & 0 & 1\\
1 & 0 & 0 & 0\\
0 & 1 & 0 & 0
\end{array}\right),\quad\quad g_1=\left(\begin{array}{cccc}
1 & 0 & 0 & 0\\
0& 1 & 0 & 0\\
0 & 0&-1 & 0\\
0 & 0 & 0 &-1
\end{array}\right).$$

%
%
%
%


\begin{example}
For any $k\in \R$, $$x(u,v)=(u, k\cosh v,v, k\sinh v)$$
defines a slant submanifold of type 2 in $(\R^4,J_1,g_1)$, with $P^2=\fra{1}{k^2+1}Id$.
\end{example}

\begin{example}\label{ocho}
For any $k \in \R$, $k\neq 0$,
$$x(u,v)= (e^{k u}\cos u \cosh v,e^{k u} \sin u \cosh v,
e^{k u} \cos u \sinh v, e^{k u} \sin u \sinh v)$$
defines a slant submanifold of type 2 in $(\R^4,J_1,g_1)$, with $P^2=\fra{k^2}{1+k^2}Id$.
\end{example}

\begin{example}\label{nueve}
For any $k \in \R$,
$$x(u,v)= (e^{k u}\cosh u \cos v,e^{k u} \sinh u \cos v,
e^{k u} \cosh u \sin v, e^{k u} \sinh u \sin v)$$
defines a totally real submanifold in $(\R^4,J,g)$.
\end{example}


Again, we can present more easy examples:
\begin{example}\label{comancheria}
For any $a,b$ with $a^2-b^2\neq 1$, $$x(u,v)=(u, av,bv, v),$$ defines a slant submanifold in $(\R^4,J_1,g_1)$, with $P^2=\fra{b^2}{1-a^2+b^2}Id$, and it is
\begin{itemize}
\item[1)] of type 1 if $a^2-b^2<1$ and $a^2>1$,
\item[2)] of type 2 if $a^2-b^2<1$ and $a^2<1$,
\item[3)] space-like of type 3 if $a^2-b^2>1$.
\end{itemize}
Moreover, it also defines a slant submanifold in $(\R^4,J,g)$, with $P^2=\fra{a^2}{1+a^2-b^2}Id$, and it is
\begin{itemize}
\item[1)] of type 1 if $b^2-a^2<1$ and $a^2>1$,
\item[2)] of type 2 if $b^2-a^2<1$ and $a^2<1$,
\item[3)] space-like of type 3 if $b^2-a^2>1$.
\end{itemize}
\end{example}

\

We can present many different examples of type 3, considering $\R^4$ with another different para Kaehler structure with metric:
$$ g_2=\left(\begin{array}{cccc}
-1 & 0 & 0 & 0\\
0& 1 & 0 & 0\\
0 & 0&1 & 0\\
0 & 0 & 0 &-1
\end{array}\right).$$

\begin{example}\label{fellini7.5}
For any $k>1$,
$$x(u,v)=(u, k\cosh v,v, k\sinh v),$$  defines a time-like slant submanifold of type 3 in $(\R^4,J_1,g_2)$, with $P^2=\fra{1}{1-k^2}Id$.
\end{example}

\begin{example}
For any $k>1$ ($k<1$),
$$x(u,v)= (e^{k u}\cosh u \cosh v,e^{k u} \sinh u \cosh v,
e^{k u} \cosh u \sinh v, e^{k u} \sinh u \sinh v)$$
defines a slant submanifold of type 1 (space-like type 3) in $(\R^4,J,g_2)$, with $P^2=\fra{k^2}{k^2-1}Id$.
\end{example}

But all the examples of type 3 slant submanifolds we have given are space-like or time-like. To offer a neutral one, we must look at least for a 4 dimensional submanifold. Combining Examples \ref{moonlight} and \ref{comancheria},  we obtain:
\begin{example}
For any $a,b$ with $a^2+b^2< 1$, $$x(u,v,z,t)=(u, av,bv, v,z, a t,\sqrt{2-b^2} t,t),$$ defines a neutral slant submanifold of type 3 in $(\R^8,J,g)$, with $P^2=\fra{a^2}{-1+a^2+b^2}Id$.
\end{example}

And with a proper election of the coefficients like in Example \ref{lalaland}, we have:

\begin{example}
For any $\theta>0$, $$x(u,v,z,t)=(u, \sinh^2\theta v,\sqrt{1-\cosh^2\theta}v, v,z, \sinh^2\theta t,\sqrt{1+\cosh^2\theta} t,t),$$ defines a neutral slant submanifold of type 3 in $(\R^8,J,g)$, with slant angle $\theta$.
\end{example}
\section{Slant curves in two dimensional Norden manifolds}
The examples given at $\R^2$ with $$ J=\left(\begin{array}{cc}
0 & 1\\
-1& 0
\end{array}\right)\quad\mbox{and}\quad  g=\left(\begin{array}{cc}
1 & 0\\
0& -1
\end{array}\right),$$
where type 3 slant submanifolds and straight lines. 

First of all, only type 3 slant curves can exist. Let $\gamma$ be a curve in a Norden manifold, and denote by ${\bf t}$ the tangent vector. Then, if $\gamma$ is a slant submanifold, $P{\bf t}$ is collineal with ${\bf t}$, therfore both of them have the same causal character and $\gamma$ is a type 3 slant submanifold.

Moreover, if the enviroment is a two dimensional manifold is Kaehler Norden we have proved the only slant curves are geodesic or invariant.

\begin{thm}
Let $\gamma$ be a smooth curve parametrized by its arc length of a 2-dimensional Kaehler Norden manifold, $\widetilde M^2$. If $\gamma$ is a slant submanifold, then it is invariant or a geodesic.
\end{thm}

\begin{proof}

The proof is similar to Theorem 2 at \cite{acslantproduct} and hence it could be ommitted.
\end{proof}



\begin{thebibliography}{99}
\bibitem{acslantpara}
{\sc P. Alegre and A. Carriazo.}
\newblock{Slant Submanifolds of Para-Hermitian Manifolds.}
\newblock{\em Mediterr. J. Math.} (2017) 14, 1-14. DOI 10.1007/s00009-017-1018-3

\bibitem{acslantproduct}
{\sc P. Alegre and A. Carriazo.}
\newblock{Curves as slant submanifolds of an almost product Riemannian manifold.}
\newblock{\em Turk. J. Math.} (2024) {\bf 48} , 701-712. https://doi.org/10.55730/
1300-0098.3535






\bibitem{borowiec}
{\sc A. Borowiec and M. Francaviglia and I. Volvovich.} 
\newblock{Anti-K\"{a}hlerian manifolds.} \emph{Differential
Geometry and its Applications}, {\bf 12} (2000), 281–289.



\bibitem{slantchen}
{\sc B.-Y. Chen.}
\newblock{Slant inmersions.}
\newblock{\em Bull. Austral. Math. Soc.} {\bf 41} (1990), 135-147.


\bibitem{Chen2}
{\sc B.-Y. Chen and O. Garay.} Classification of quasi-minimal surfaces with parallel mean
curvature vector in pseudo-Euclidean 4-space $\mathbb{E}^4_2$. \emph{Results Math.} \textbf{55} (2009), 23-38.

\bibitem{Chen3}
{\sc B.-Y. Chen and I. Mihai.} Classification of quasi-minimal slant surfaces in Lorentzian
complex space forms. \emph{Acta Math. Hungar.} \textbf{122} (2009), 307-328.

\bibitem{etayo}
{\sc F. Etayo and R. Santamar\'ia.} $(J^2=\pm I)$-metric manifolds. \emph{Publ. Math. Debrecen.} {\bf 57/3-4} (2000), 435-444.

\bibitem{ganchev}
{\sc G.T. Ganchev and A.V. Borisov.} Note on the almost complex manifolds with a Norden metric. \emph{C.R. Acad. Bulgarie.} {\bf 39} (1986), 31-34.

\bibitem{ganchev2}
{\sc G. Ganchev and K. Gribachev and V. Mihova.} B-connections and their conformal invariants on conformally
Kaehler manifolds with B-metric. \emph{Publ. de L’Inst. Math.} {\bf 42} (56) (1985), 107–121.

\bibitem{norden}
{\sc A.P. Norden.} On a class of four-dimensional $A$-spaces. \emph{Izv. Vys\v{s}. U\v{c}ebn. Zava ested. Matematika.} {\bf 4} (1960), 145-157.

\bibitem{chr}
{\sc R. Castro, L.M. Hervella and E.G. Rio.} Some examples of almost complex manifolds with Norden metric. \emph{Riv. Mat. Univ. Parma (4).} {\bf 15} (1989), 133-141.

\bibitem{gupta}
{\sc B.K. Gupta and B.B: Chaturvedi.} On a slant submanifold of a Kaehler-Norden manifold. \emph{J. Tensor Soc.} {\bf 16} (2022), 71-80.

\bibitem{ay}
{\sc V. Ayhan and S.P. Perkta\c{s}.} Slant submanifolds of almost poly-Norden Riemannian manifolds. \emph{Int. J. Maps Math.} {\bf 6}(2023), no. 1, 22-36.







\bibitem{sahin}
{\sc B. Sahin.}
\newblock{Slant submanifolds of an almost product Riemannian manifold.}
\newblock{\em J. Korean Math. Soc.} {\bf 43} (2006), No. 4, 717-732.

\end{thebibliography}
\end{document}